\documentclass[a4paper,11pt]{article}

\usepackage[left=3cm,top=3cm,right=3cm, bottom=4cm]{geometry}
\usepackage{booktabs}
\usepackage{tabularx}
\usepackage{amsthm}
\usepackage{amsmath}
\usepackage{amssymb}
\usepackage{subfig}
\usepackage{graphicx}
\usepackage{verbatim}
\usepackage{appendix}
\usepackage{listings}
\usepackage{enumerate}
\usepackage{natbib}
\usepackage{cases}
\usepackage{bbm}
\usepackage{color}

\newcommand{\E}{\mathbb{E}}
\newcommand{\nn}{\nonumber}
\newcommand{\p}{\mbox{P}}

\title{Fluid Approximation of a Call Center Model with Redials and Reconnects}
\author{Sihan Ding, Maria Remerova, Rob van der Mei, Bert Zwart}

\begin{document}
\maketitle

\begin{abstract}
In many call centers, callers may call multiple times. Some of the calls are re-attempts after abandonments (redials), and some are re-attempts after connected calls (reconnects). The combination of redials and reconnects has not been considered when making staffing decisions, while ignoring them will inevitably lead to under- or overestimation of call volumes, which results in improper
and hence costly staffing decisions.
Motivated by this, in this paper we study call centers where customers can abandon, and abandoned customers may redial, and when a customer finishes his conversation with an agent, he may reconnect. We use a fluid model to derive first order approximations for the number of customers in the redial and reconnect orbits in the heavy traffic. We show that the fluid limit of such a model is the unique solution to a system of three differential equations. Furthermore, we use the fluid limit to calculate the expected total arrival rate, which is then given as an input to the Erlang A model for the purpose of calculating service levels and abandonment rates. The performance of such a procedure is validated in the case of single intervals as well as multiple intervals with changing parameters.
\end{abstract}

\section{Introduction}
Nowadays, call centers are important means of communication with the customer. Therefore, the response-time performance of call centers is crucial for the
customer satisfaction. For call center managers, making the right staffing decisions (i.e., decide on the right number of agents) is essential to the costs and the performances of call centers. Various models have been developed in order to decide on the right number of agents. One of the most widely used models is the Erlang C model and there is a lot of literature on it (see \citet{gans2002telephone} and the references therein). A staffing rule called the square-root staffing is proposed by \citet{halfin1981heavy}. \citet{garnett2002designing} show that the square-root staffing rule remains valid for the Erlang A model. However, both the Erlang C formula and the square-root staffing formula ignore customer redial (a re-attempt after an abandoned call) behaviors in call centers, while this behavior is quite significant (see \citet{gans2002telephone} and reference therein). \citet{aguir2008interaction} discover that ignoring redials can lead to under-staffing or over-staffing, depending on the forecasting assumption being made. \citet{sze1984} studies a model where abandonments and redials are included, focusing on the traffic loaded systems.

Besides redials, there also exists another important feature, which is called reconnect (a re-attempt after a connected call). The reconnect customer behavior is first mentioned in \citet{gans2002telephone}, where the reconnect is defined as a revisit. In \citet{ding2013}, we use real call center data to show that an inbound call can either be a fresh call (a initial attempt), a redial or a reconnect. Also, as argued in \citet{ding2013}, redials and reconnects should be considered and modeled, since ignoring them can lead to significantly inaccurate estimations of the total inbound volume. As a consequence, neglecting the impact of redials and reconnects will lead to either overstaffing or understaffing. In case of overstaffing the performance of the call center will be good,
but at unnecessarily high costs. In case of understaffing, the performance of the call center will be degraded, which may lead to customer dissatisfaction and possibly customer churn. Despite the economic relevance of including both features in staffing models,
to the best of the authors' knowledge no papers have appeared on staffing of call centers where {\it both} redials and reconnects are included. This paper aims to fill this gap, that is, we investigate the staffing problem in call centers with the features of both redials and reconnects. We focus on the case of large call centers that operate under heavy load.

Intuitively, when the system is heavily loaded, it would lead to bad service levels (SL). However, for large call centers, especially during the busy hours when the inbound volume is quite large, it is possible that the target SL can be met even in heavy traffic. Further discussions of this effect can be found in \citet{garnett2002designing} and \citet{borst2004}.

In this paper, we aim to answer the following question: ``In large call centers, what are the SL and the abandonment percentage (AP) if both redialling and reconnection of customers are taken into account?" To this end, one must first estimate the total number of arrivals into the call center. This is not trivial, since the number of total arrivals depends on the number of agents (see \citet{ding2013}). This dependency effect between the total number of arrivals and the number of agents becomes more complicated in real life, due to the fact that the fresh arrival rate and the number of agents are often time-dependent. If the number of arrivals cannot be determined, it is impossible to calculate the SL. Therefore, in this paper, we take a two-step approach to calculate the SL and AP. First, we numerically calculate the expected total arrival rate at any instant time by using a fluid limit approximation. We also show that the fluid limit of this model is a unique solution to a system of three deterministic differential equations. In the second step, under the assumption of the total arrival process being Poisson, we apply the Erlang A formula to obtain the SL and the AP. This approximation turns out to be quite accurate. In this paper, we consider only the expected SL and AP, for discussion about the SL variability, we refer to the work by \citet{roubos2012service}.

Fluid models for call centers have been extensively studied. \citet{whitt2006fluid} gave an intuitive explanation of the fluid model. He develops a deterministic fluid limit which they use to provide first-order performance descriptions for the $G/GI/s+GI$ queueing model under heavy traffic, where the second $GI$ stands for the i.i.d. patience distribution. In \citet{whitt2006fluid}, the redial behavior is not considered, though. The existence and uniqueness of the fluid limit are given as conjectures. \citet{mandelbaum2002queue} use the fluid and diffusion approximation for the multi-server system with abandonments and redials. He obtains first order approximations of queue length and expected waiting time as well as their confidence bounds. In \citet{mandelbaum1999time}, the authors use the fluid and the diffusion approximation for time varying multiserver queue with abandonment and retrials. They show that the fluid and the diffusion approximation can both be obtained by solving sets of non-linear differential equations, where the diffusion process can provide the confidence bounds for the fluid approximation. The work by \citet{mandelbaum1998strong} gives more general theoretical results for the fluid and diffusion approximation for Markovian service networks. \citet{aguir2004impact} extend the model by allowing customer balking behavior, but no formal proof of the fluid limit is given. Besides the application in call center staffing problems, fluid models have also been applied in delay announcement of customers in call centers (see \citet{ibrahim2009real, ibrahim2011wait}).

The rest of the paper is structured as follows. In section~\ref{sec:modeldescription}, we describe the queueing model with the features of the redial and reconnect. In section~\ref{sec:fluidmodel}, we propose a fluid model, which is a deterministic analogue of the stochastic model. We prove that the original stochastic model converges to the fluid model under a proper scaling. We numerically compute the fluid model, and simulate the original model, and compare them in the case of a single interval and multiple intervals, where the parameters are changing per interval and would remain piece-wise constants within each interval. The Erlang A formula is then used to approximate the waiting time distributions in section~\ref{erlanga}.

\section{Model description} \label{sec:modeldescription}
\begin{figure}[ht]
  \begin{center}
    \includegraphics[width=0.8\textwidth]{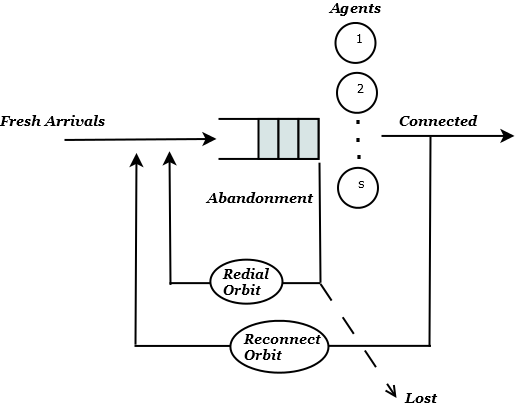}
  \caption{Call diagram}
  \label{fig:diagram}
  \end{center}
\end{figure}

Consider the queueing model illustrated in Figure~\ref{fig:diagram}. We assume that calls arrive according to a Poisson process. We refer to these calls as {\it fresh calls}. There are $s$ agents who handle inbound calls. An arriving call is handled by an available agent, if there is any; otherwise, he waits in an infinite buffer queue. The calls are handled in the order of arrival. After an exponentially distributed amount of time $\Psi$, a waiting customer who did not get connected to an agent will lose his patience and abandon. We assume $\E \Psi =1/ \theta < \infty$, where $\theta$ is the abandonment rate. With probability $p$, an abandoned customer will enter the redial orbit, and he will redial after an exponentially distributed amount of time $\Gamma_{RD}$, with $\E \Gamma_{RD} = \delta_{RD} < \infty$. We refer to these calls as {\it redials}. With probability $1-p$, this customer will not call back, and this call is considered as a ``lost'' call. We assume that the service time $B$ of a customer has an exponential distribution with mean $\E B = 1/\mu < \infty$. After the call has been served, this customer will enter the reconnect orbit with probability $q$, and he will reconnect after an exponentially distributed time $\Gamma_{RC}$, with $\E \Gamma_{RC} = \delta_{RC} < \infty$. We refer to such calls as {\it reconnects}. We assume that $p$ and $q$ do not depend on customers' experiences in the system. These experiences include holding time, waiting time and the number of times that customers have already called. We use this queueing model to represent the situation of a single-skill call center. In this paper, we consider independent service times; for the study of dependent service times, please see \citet{pang2012impact}.

\section{Fluid limit approximations} \label{sec:fluidmodel}
In this section, we first show that the problem of calculating the expected total arrival rate drills down to the problem of calculating $\E Z_{Q}(t), \E Z_{RD}(t)$ and $\E Z_{RC}(t)$, where $Z_{Q}(t)$ is the number of customers in the queue plus the number of customers in service at time $t$, $Z_{RD}(t)$ is the number of customers in the redial orbit at time $t$, and $Z_{RC}(t)$ is the number of customers in the reconnect orbit at time $t$.  Because an arrival can be a fresh arrival, a redial or a reconnect, the following equation must hold
\begin{align} \label{eq:totalarrival}
  \E \Lambda(t) &= \lambda (t) + \E \lambda_{RD}(t) + \E \lambda_{RC}(t) \nn \\
                &= \lambda (t) + \delta_{RD} \E Z_{RD}(t)  + \delta_{RC} \E Z_{RC}(t),
\end{align}
where $\Lambda(t)$ stands for the total arrival rate at time $t$, which is a stochastic process, $\lambda(t)$ stands for the fresh arrival rate at time $t$, $\lambda_{RD}(t)$ and $\lambda_{RC}(t)$ stand for the arrival rate due to redials and reconnects at time $t$, respectively. Therefore, once $\E Z_{Q}(t), \E Z_{RD}(t)$ and $\E Z_{RC}(t)$ are known, $\E \Lambda(t)$ can be obtained by Equation \eqref{eq:totalarrival}. Note that $Z_{Q}(t)$ does not appear in Equation \eqref{eq:totalarrival}, but we will see later that $Z_{RD}(t)$ and $Z_{RC}(t)$ would depend on $Z_{Q}(t)$.

In fact, the stochastic process $\{\mathbf{Z}(t), t \geq 0 \}$, which is defined by
\begin{align} \label{def1}
\mathbf{Z}(t) := \left( Z_{Q}(t), Z_{RD}(t), Z_{RC}(t) \right)^T,
\end{align}
is a 3-dimensional Markov process, because the inter-arrival time, service duration and other durations are assumed to be exponentially distributed. The state space of this Markov process is $\mathbb{Z}_+^3$. To save space, we will not show the transition diagram here. Since it is a Markov process, we can truncate the system at certain large state, and numerically obtain the steady state distribution of $\mathbf{Z}(t)$ by solving global balance equations. Theoretically, by truncating at some large states, this method offers almost exact results. However, for the model we consider, it is very difficult to formulate and solve the global balance equations, and their solution offers no insight about the system. Therefore, for the convenience of practical usage, we will not consider solving this Markov process, but other approximation methods.

\subsection{Fluid limit}
In this subsection, we present the fluid model, which we show to arise as the limit under a proper scaling of the stochastic model in Figure \ref{fig:diagram}.

Consider a single interval with the fresh arrival rate remaining constant during this interval (e.g., $\lambda(t) = \lambda$, $t \geq 0$). The following flow conservation equations hold for this stochastic model:
\begin{align}
    Z_{Q} \left( t \right) & = Z_Q \left( 0 \right) + \Pi_{\lambda}\left( t \right) + D_{RD} \left( t \right) + D_{RC}\left( t \right) - D_s \left( t \right) - D_{a}\left( t \right),  \label{eq:flowq}\\
    Z_{RD}\left( t \right) & = Z_{RD}\left( 0 \right) + \sum_{j=1}^{D_a\left( t \right)} B_j\left( p \right) - D_{RD}\left( t \right), \label{eq:flowrd}\\
    Z_{RC}\left( t \right) & = Z_{RC} \left( 0 \right) + \sum_{j=1}^{D_s\left( t \right)} B_j\left( q \right) - D_{RC}\left( t \right), \label{eq:flowrc}
\end{align}
where $\Pi_{\lambda}\left( t\right)$ is the number of fresh arrivals during time interval $[0, t)$, and $\Pi_{\lambda}(\cdot)$ is a Poisson process of rate~$\lambda$. In addition, $D_{RD} \left( t \right), D_{RC} \left( t \right), D_s \left( t \right),  D_{a} \left( t \right)$ are the number of redials during $[0, t)$, number of reconnects during $[0, t)$, number of served customers during $[0, t)$ and number of abandoned customers during time $[0, t)$, respectively. $B_j\left( p \right)$ is a Bernoulli random variable with success probability $p$, $j = 1, 2, \ldots, D_a\left( t \right)$. $B_j\left( p \right) = 1$, if the $j$th abandoned customer enter the redial orbit; $B_j\left( p \right) = 0$, otherwise. Therefore, for given $D_a\left( t \right)$, $\sum_{j=1}^{D_a\left( t \right)} B_j\left( p \right) \sim \operatorname{Bin} (D_a\left( t \right), p)$. By the same argument, we have $\sum_{j=1}^{D_s\left( t \right)} B_j\left( q \right) \sim \operatorname{Bin} (D_s\left( t \right), q)$.

Let $\Pi_{i}(\cdot)$, $i = 1,2, 3, 4$, be independent Poisson processes of rate 1, then we claim the following
\begin{align*}
    D_s(t) & = \Pi_{1}\left(\int_{0}^{t}\mu \min \{s, Z_{Q}(u)\} du\right), \nn \\
    D_a(t) & =  \Pi_{2}\left(\int_{0}^{t}\theta \left( Z_{Q}\left( u \right) - s\right)^+ du\right), \nn \\
  D_{RD} \left( t \right) & =\Pi_{3}\left(\int_{0}^{t}\delta_{RD} Z_{RD}\left( u \right) du\right), \nn \\
  D_{RC} \left( t \right) &=\Pi_{4}\left(\int_{0}^{t}\delta_{RC} Z_{RC}\left( u \right) du\right). \nn
\end{align*}
Rigorous proof of these four statements can be given along the lines of \citet{pang2007martingale}, see Lemma 2.1.

To introduce the fluid limit, we consider a sequence of models as in Figure \ref{fig:diagram} such that, in the $n$-th model, the fresh arrival rate is $\lambda n$ and the number of servers is $ns$. We add the superscript $``\left( n\right)"$ to all notations in the $n$-th model. Similarly to \eqref{eq:flowq}-\eqref{eq:flowrc}, we then have for the $n$-th model:
\begin{align}
     Z_{Q}^{\left( n\right)} \left( t \right) & = Z_Q^{\left( n \right) } \left( 0 \right) + \Pi_{\lambda n}^{\left( n\right)} \left( t \right) + D_{RD}^{ \left( n \right)}  \left( t \right) + D_{RC}^{ \left( n \right)} \left( t \right) - D_s^{ \left( n \right)}  \left( t \right) - D_{a}^{ \left( n \right)}  \left( t \right),  \label{eq:fluidqn}\\
     Z_{RD}^{ \left( n \right)} \left( t \right) & = Z_{RD}^{ \left( n \right)} \left( 0 \right) + \sum_{j=1}^{D_a^{\left( n \right)}\left( t \right)} B_j\left( p \right) - D_{RD}^{\left( n \right)}\left( t \right), \label{eq:fluidrdn}\\
     Z_{RC}^{ \left( n \right)} \left( t \right) & = Z_{RC}^{ \left( n \right)} \left( 0 \right) + \sum_{j=1}^{D_s^{\left( n \right)}\left( t \right)} B_j\left( q \right) - D_{RC}^{\left( n \right)}\left( t \right). \label{eq:fluidrcn}
\end{align}

Now we define the fluid scaled process \begin{align*}
\bar{\mathbf{Z}}^{(n)}(t):= \left( \bar{Z}_{Q}^{\left( n\right)} \left( t \right), \bar{Z}_{RD}^{\left( n\right)} \left( t \right), \bar{Z}_{RC}^{ \left( n \right)} \left( t \right) \right)^T,
\end{align*}
 where
\begin{align*}
\bar{Z}_{Q}^{\left( n\right)} \left( t \right) := \frac{Z_{Q}^{\left( n\right)} \left( t \right)}{ n}, \hspace{3mm} \bar{Z}_{RD}^{\left( n\right)} \left( t \right) := \frac{Z_{RD}^{ \left( n \right)} \left( t \right)}{n}, \hspace{3mm} \bar{Z}_{RC}^{ \left( n \right)} \left( t \right) := \frac{Z_{RC}^{ \left( n \right)} \left( t \right)}{n}.
\end{align*}

Let $D([0, \infty), \mathbb{R}^3)$ be the space of right continuous functions with left limits in $ \mathbb{R}^3$ having the domain $[0, \infty)$. We endow $D([0, \infty), \mathbb{R}^3)$ with the usual Skorokhod $J_1$ topology. Suppose $\{X^{(n)}\}_{n=1}^{\infty}$ is a sequence of stochastic processes, then notation $X^{(n)} \xrightarrow{d} x$ means that $X^{(n)}$ converge weakly to stochastic process $x$.
\newtheorem{def1}{Definition}
\begin{def1} \label{def:def1}
If there exists a limit in distribution for the scaled process $\{\bar{\mathbf{Z}}^{(n)}(\cdot)\}_{n=1}^{\infty}$, i.e. $\bar{\mathbf{Z}}^{(n)}(\cdot) \xrightarrow{d} \mathbf{z}(\cdot)$, then $\mathbf{z}(\cdot)$ is called the fluid limit of the original stochastic model.
\end{def1}

\subsubsection{Fluid limit for a single interval}
To obtain the fluid limit of the system (i.e., a sequence of stochastic processes specified by Equations \eqref{eq:fluidqn}-\eqref{eq:fluidrcn}) for a single interval, we divide both sides of Equations \eqref{eq:fluidqn}-\eqref{eq:fluidrcn} by $n$, then let $n \rightarrow \infty$.

\newtheorem{theo2}{Lemma}
\begin{theo2} \label{lemma}
The sequence of scaled processes $\displaystyle \{\mathbf{\bar{Z}}^{\left( n \right)}( \cdot)\}_{n=1}^{\infty}$ is relatively compact and all weak limits are a.s. continuous.
\end{theo2}
\begin{proof}
See Appendix A.2.
\end{proof}

\newtheorem{theo1}{Theorem}
\begin{theo1} \label{theo:theo1}
If for given deterministic values $\left( z_{Q}(0), z_{RD}(0), z_{RC}(0) \right) $, we assume \\
$\left( \bar{Z}_Q^{\left( n\right)}(0), \bar{Z}_{RD}^{\left( n\right)}(0), \bar{Z}_{RC}^{\left( n\right)}(0) \right) \xrightarrow{d} \left( z_{Q}(0), z_{RD}(0), z_{RC}(0) \right)$ as $n \rightarrow \infty$, then the fluid limit of the original stochastic model is the unique solution to the following system of equations
\begin{align}
  z_{Q}(t) &= z_{Q} (0) + \lambda t + \delta_{RD} \int_{0}^{t} z_{RD}(u) du +  \delta_{RC} \int_{0}^{t} z_{RC}(u) du  \nn \\
  & - \mu \int_{0}^{t} \min \{ s, z_Q(u)\} du - \theta \int_0^t \left( z_Q(u) - s \right)^+ du,  \label{eq:fluidlimitq}\\
  z_{RD}(t) &= z_{RD}(0) + p\theta \int_0^t \left( z_Q(u) - s \right)^+ du - \delta_{RD} \int_{0}^{t} z_{RD}(u) du, \label{eq:fluidlimitrd}\\
  z_{RC}(t) &= z_{RC}(0) + q\mu \int_{0}^{t} \min \{ s, z_Q(u)\} du - \delta_{RC} \int_{0}^{t} z_{RC}(u) du \label{eq:fluidlimitrc}.
\end{align}
\end{theo1}
\begin{proof}
See Appendix A.3.
\end{proof}

We could not obtain analytic expressions of $z_{Q}(t), z_{RD}(t)$ and $z_{RC}(t)$ from Equations \eqref{eq:fluidlimitq}-\eqref{eq:fluidlimitrc}. However, we can solve them numerically by the following iterative procedure. First, Equations \eqref{eq:fluidlimitq}-\eqref{eq:fluidlimitrc} can be rewritten as
\begin{align*}
\mathbf{z}(t) = \Phi(\mathbf{z}(t)).
\end{align*}
We let $\mathbf{z}^{(0)}(0)= 0$, and then calculate $\mathbf{z}^{(k+1)} = \Phi(\mathbf{z}^{(k)})$, $k = 0, 1, \ldots$, until the difference between $\mathbf{z}^{(k+1)}$ and $\mathbf{z}^{(k)}$ is small enough.

\subsubsection{Fluid limit for multiple intervals}
We have just shown the fluid limit for a single interval, where the parameters $\lambda$ and $s$ remain the same within the interval. However, in real call centers, parameters can vary during the day, especially the arrival rate $\lambda(t)$. As shown by \citet{shen2008} and \citet{ibrahim2013}, call volumes normally follow certain intraday patterns. Observing the intraday arrival pattern from the historical data set, managers would schedule different number of agents for each interval to meet the SL. Therefore, we now show the fluid limit for multiple intervals, where $\lambda$ and $s$ vary vary from interval to interval. We assume that other parameters remain constant.

We divide a day into $m$ intervals. Each interval starts at $t_{i-1}$ and ends at $t_i$, $i = 1, 2, \ldots, m$. The fresh arrival rate of interval $i$ is denoted by $\lambda_i$, and the number of agents in interval $i$ is denoted by $s_i$, $i = 1, 2, \ldots, m$. For the $i$th interval, i.e., $t_{i-1} \leq t < t_i$, the fluid limit then becomes
\begin{align}
  z_{Q}(t) &= z_{Q} (t_{i-1}) + \lambda_i \left( t-t_{i-1}\right) + \delta_{RD} \int_{t_{i-1}}^{t} z_{RD}(u) du +  \delta_{RC} \int_{t_{i-1}}^{t} z_{RC}(u) du  \nn \\
  & - \mu \int_{t_{i-1}}^{t} \min \{ s_i, z_Q(u)\} du - \theta \int_{t_{i-1}}^{t} \left( z_Q(u) - s_i \right)^+ du,  \label{eq:fluidlimitq2}\\
  z_{RD}(t) &= z_{RD}(t_{i-1}) + p\theta \int_{t_{i-1}}^t \left( z_Q(u) - s_i \right)^+ du - \delta_{RD} \int_{t_{i-1}}^{t} z_{RD}(u) du, \label{eq:fluidlimitrd2}\\
  z_{RC}(t) &= z_{RC}(t_{i-1}) + q\mu \int_{t_{i-1}}^{t} \min \{ s_i, z_Q(u)\} du - \delta_{RC} \int_{t_{i-1}}^{t} z_{RC}(u) du \label{eq:fluidlimitrc2}.
\end{align}

Numerically solving Equations \eqref{eq:fluidlimitq2}-\eqref{eq:fluidlimitrc2} is similar to solving Equations \eqref{eq:fluidlimitq}-\eqref{eq:fluidlimitrc}, thus, we do not elaborate on the procedure here.

In reality, parameters such as $\mu, \theta, \delta_{RD}$ and $\delta_{RC}$ can also be time-dependent, and vary per interval. For example, $\delta_{RD}$ may be bigger in the late afternoon than in the morning, since abandoned customers want to have responses by the end of the day. It is possible to extend the model in Equations \eqref{eq:fluidlimitq2}-\eqref{eq:fluidlimitrc2} to adapt such situation by simply replacing the parameters. In this paper, for the simplicity of validation, we will not consider such cases.

\subsection{Model under stationarity}
We have just shown that one can numerically solve differential equations \eqref{eq:fluidlimitq}-\eqref{eq:fluidlimitrc} to obtain the fluid limit $\mathbf{z}(t)$. We now derive the fluid limit in stationarity, i.e., we develop conditions under which $\mathbf{z}(t)$ is constant.

By differentiating Equations \eqref{eq:fluidlimitq}-\eqref{eq:fluidlimitrc} and taking into account that $\displaystyle{\frac{d}{dt}\mathbf{z}(t) = 0}$ for a constant solution, we obtain
\begin{align}
  0 &=  \lambda + \delta_{RD} z_{RD}(\infty) +  \delta_{RC}  z_{RC}(\infty) - \mu \min \{ s, z_Q(\infty)\} - \theta \left( z_Q(\infty) - s \right)^+, \label{eq:stationaryq} \\
  0 &=  p\theta \left( z_Q(\infty) - s \right)^+ - \delta_{RD} z_{RD}(\infty), \label{eq:stationaryrd}\\
  0 &=  q\mu \min \{ s, z_Q(\infty)\} - \delta_{RC} z_{RC}(\infty) \label{eq:stationaryrc},
\end{align}
where $\displaystyle z_{Q}(\infty):=\lim_{t \rightarrow \infty}z_{Q}\left( t\right), z_{RD}(\infty):=\lim_{t \rightarrow \infty}z_{RD}\left( t\right), z_{RC}(\infty):=\lim_{t \rightarrow \infty}z_{RC}\left( t\right)$.

Equations \eqref{eq:stationaryq}-\eqref{eq:stationaryrc} can be easily solved with respect to $z_{Q}(\infty), z_{RD}(\infty)$ and $z_{RC}(\infty)$, yielding
\begin{align}
  z_Q(\infty) &= \begin{cases} \label{eq:zqcase}
        \displaystyle{\frac{\lambda}{\left( 1 - q \right) \mu}},  \hspace{10mm} \text{if }\rho < \left( 1 - q \right) \\
        \displaystyle{\frac{\lambda + q\mu s -\mu s}{\theta \left( 1 - p \right)} + s}, \hspace{10mm} \text{if } \rho \geq \left( 1 - q \right) \\
        \end{cases}\\
  z_{RD}(\infty) &= \begin{cases} \label{eq:zrdcase}
        \displaystyle{\frac{p \theta \left( z_Q(\infty)-s \right) }{\delta_{RD}}},  \hspace{10mm} \text{if }\rho < \left( 1 - q \right) \\
        0, \hspace{10mm} \text{if } \rho \geq \left( 1 - q \right)
        \end{cases} \\
  z_{RC}(\infty) &= \begin{cases} \label{eq:zrccase}
        \displaystyle{\frac{q \mu z_{Q}(\infty)}{\delta_{RC}}}, \hspace{10mm} \text{if }\rho < \left( 1 - q \right) \\
        \displaystyle{\frac{q \mu s}{\delta_{RC}}}. \hspace{10mm} \text{if }\rho \geq \left( 1 - q \right) \\
        \end{cases}
\end{align}

The results above would offer some insights. $\rho: = \frac{\lambda}{s \mu}$ is the load of the system due to the fresh arrivals. Since $\frac{1}{1-q}$ portion of $\rho$ will reconnect, the total load would be $\hat{\rho}:= \frac{\lambda}{(1-q)s\mu}$, if there were no redials. One can notice that in expressions \eqref{eq:zqcase}-\eqref{eq:zrccase}, the value of $\hat{\rho}$ determines whether the fluid model is in heavy traffic or not. Thus, from now on, we use $\hat{\rho}$ to denote the actual load of the system instead of $\rho$. In the case of $\hat{\rho} < 1$, since the fluid limit is deterministic, $z_Q(\infty) < s$, and $z_{RD}(\infty) = 0$ would hold. This means there is no abandonment at all in the fluid limit when $\hat{\rho} < 1$. In reality, due to the variability of the service duration and patience, abandonments would not be $0$ though, but very small. If $\hat{\rho} > 1$, by Equation \eqref{eq:zqcase}, $z_Q(\infty) > s$. Therefore, in this case, the fluid model indicates that there will be $(z_Q(\infty) - s)$ amount of customers waiting, each with rate $\theta$, and customers will go to the redial orbit with rate $p\theta(z_Q(\infty)-s)$.

\section{Validation of the fluid limit}
In this section, we will validate the fluid model via simulation both for a single interval and for multiple intervals. We simulate the system for $480$ minutes of time, i.e., $8$ hours, which correspond to the busy hours in some call centers. The  results obtained via the fluid limit are compared with the simulation results. Since $\mathbf{Z}(t)$ is a stochastic process, it has variability. To remove those variabilities, we do the simulation for $100$ times, and then take the average.

\subsection{Validation of a single interval}
\begin{figure}[ht]
  \begin{center}
    \includegraphics[width=0.7\textwidth]{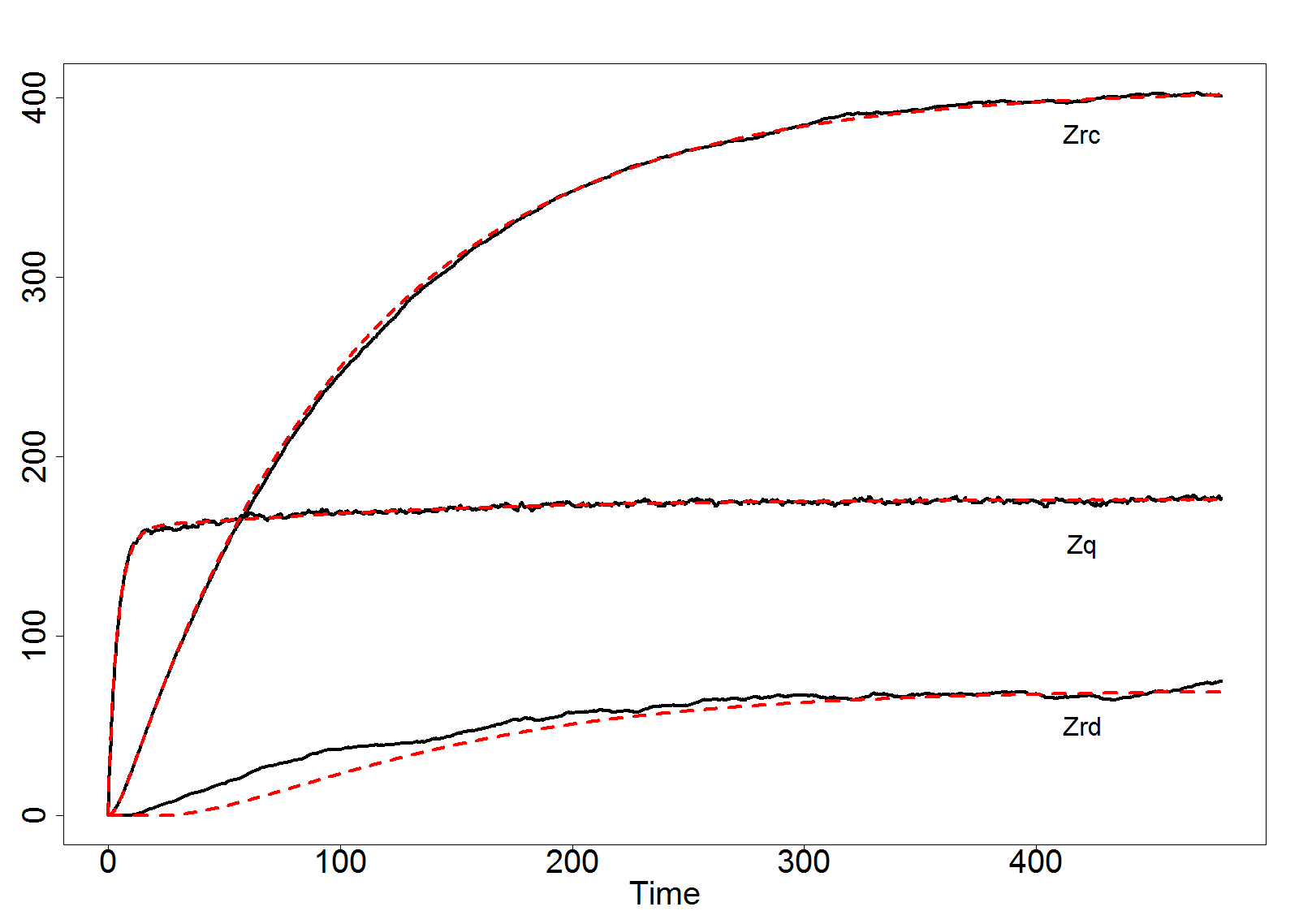}
  \caption{$z_{RC}(t), z_{Q}(t), z_{RD}(t)$ (from top to bottom) via fluid limit (red dashed curve) and $Z_{RC}(t), Z_{Q}(t), Z_{RD}(t)$ via simulation (black solid curve), $\lambda = 40, \mu=1/4, p=0.5, q=0.1, \theta = 1/2, s=148,	\delta_{RD}=0.05, \delta_{RC} = 0.01$}
  \label{fig:validation}
  \end{center}
\end{figure}

We start with the simple case of a single interval, where $\lambda(t) = \lambda$, for all $t > 0$, and we assume that $s, \mu$ as well as other parameters are constants over time. We compare $\mathbf{z}(t)$ (computed via Equations \eqref{eq:fluidlimitq}-\eqref{eq:fluidlimitrc}) with $\mathbf{Z}(t)$ (simulation results) for different values of $\hat{\rho}$. For each value of $\hat{\rho}$, $s$ changes, while $\lambda$ and $\mu$ remain the same.  One example of $\mathbf{z}(t)$ and $\mathbf{Z}(t)$, where $\hat{\rho} = 1.2$, is shown in Figure~\ref{fig:validation}.

One can see from Figure~\ref{fig:validation} that the system starts with zero customers, and as time passes by, $Z_{Q}(t), Z_{RD}(t)$ and $Z_{RC}(t)$ gradually build up and reach stationarity. Furthermore, in this parameter setting, the fluid limit offers a close approximation of the original system, especially for $Z_{RQ}(t)$ and $Z_{RC}(t)$. The errors are bigger for $Z_{RD}(t)$.

Obtaining the approximation of $\mathbf{Z}(t)$ is the intermediate step for calculating $\lambda_{RD}(t)$ and $\lambda_{RC}(t)$. Therefore, for the purpose of testing the errors of the fluid model in number of redials and reconnects, we introduce the error measurements $e_{RD}$ and $e_{RC}$, which are defined by
\begin{align*}
e_{RD} := \frac{\int_{0}^{T}|\E \lambda_{RD}(u) - \lambda_{RD}^{fl}(u)|du}{\int_{0}^{T}\E \lambda_{RD}(u)du} = \frac{\int_{0}^{T}|\E Z_{RD}(u) - z_{RD}(u)|du}{\int_{0}^{T}\E Z_{RD}(u)du},\\
e_{RC} := \frac{\int_{0}^{T}|\E \lambda_{RC}(u) - \lambda_{RC}^{fl}(u)|du}{\int_{0}^{T}\E \lambda_{RC}(u)du} =
\frac{\int_{0}^{T}|\E Z_{RC}(u) - z_{RC}(u)|du}{\int_{0}^{T}\E Z_{RC}(u)du},
\end{align*}
where $\lambda_{RD}^{fl}(t)$ and $\lambda_{RC}^{fl}(t)$ are the arrival rate due to redial and reconnect in the fluid approximation, respectively, and $T = 480$, as the same length of the simulation time. The parameters and results are shown in Table \ref{table:scenarios1}.

\begin{table}[ht]
  \begin{center}
  \begin{tabular}{@{}cccc@{}}
    \toprule
    $\hat{\rho}$ & $s$ & $e_{RD}$ & $e_{RC}$  \\
   	\midrule
  	1.01 & 176 & 92.5\% & 1.7\% \\
    1.05 & 169 & 35.7\% & 1.6\% \\
    1.1 & 162 & 10.3\% & 0.5\% \\
    1.2 & 148 & 1.9\% & 0.5\% \\
    1.3 & 137 & 1.3\% & 0.5\% \\
    1.4 & 127 & 1.4\% & 0.5\% \\
    1.5 & 119 & 1.1\% & 0.7\% \\
    \bottomrule
  \end{tabular}
  \caption{Approximation errors of different values of $\hat{\rho}$ in single intervals , $\lambda = 40, \mu=1/4, p=0.5, q=0.1, \theta = 1/2,	\delta_{RD}=0.05, \delta_{RC} = 0.01$}
  \label{table:scenarios1}
  \end{center}
\end{table}
One can see from Table \ref{table:scenarios1} that for the number of reconnects,  the fluid model offers good approximations for all scenarios. However, for the number of redials, the fluid model performs badly when $\hat{\rho} < 1.1$. In the next section, we will show that the consequences of these bad performances are not severe in terms of SL and AP. When $\hat{\rho} \geq 1.1$, the fluid model starts to get more accurate with $e_{RD} \leq 10.3\%$. In cases where $\hat{\rho} \geq 1.2$, $e_{RD}$ is less then $1.1\%$.

\subsection{Validation of multiple intervals}
\begin{figure}[ht]
  \begin{center}
    \includegraphics[width=0.7\textwidth]{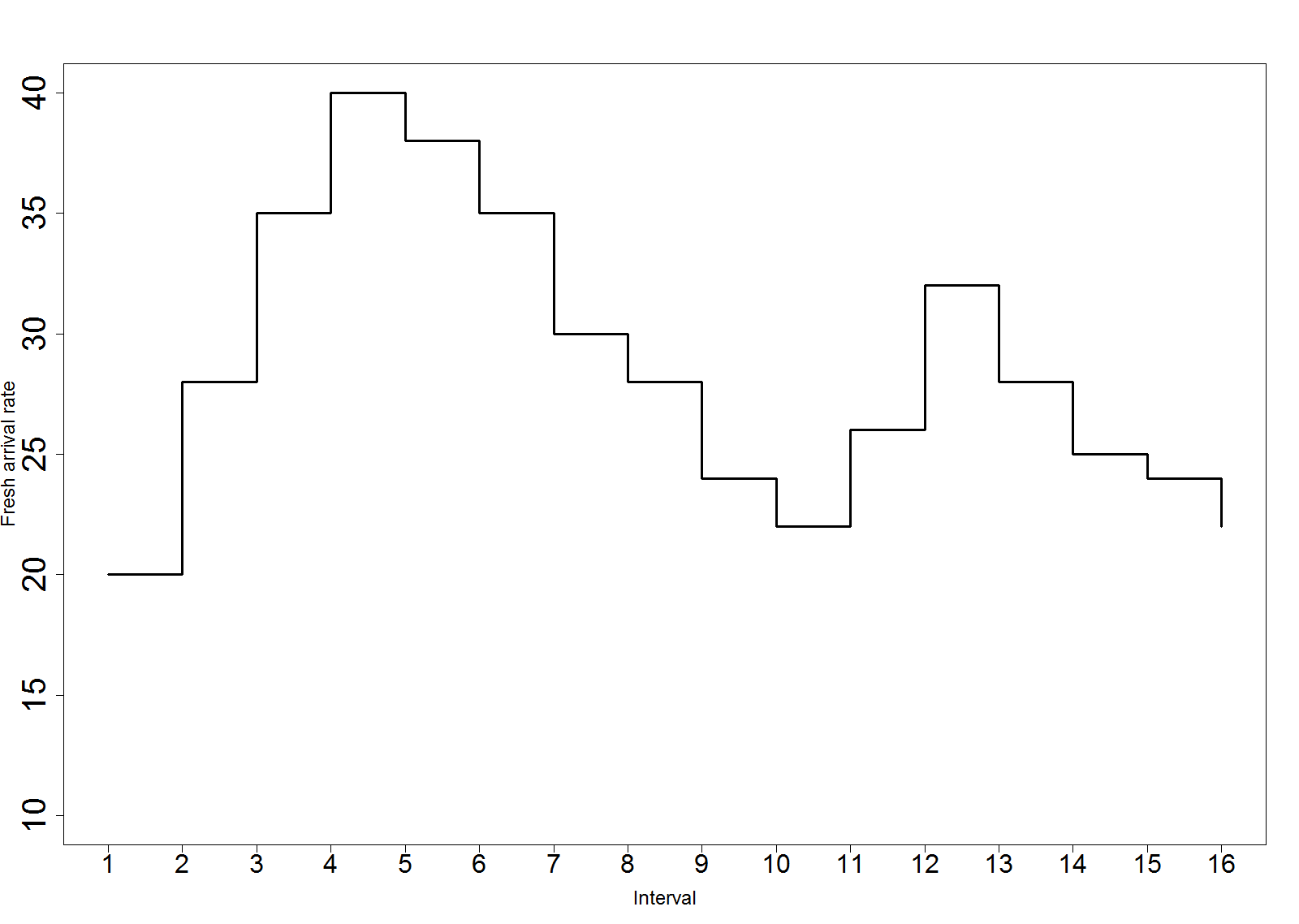}
  \caption{Fresh arrival rate per interval}
  \label{fig:arrivalpattern}
  \end{center}
\end{figure}

Similar to the validation procedure in the case of a single interval, now we validate the performance of the fluid model for multiple intervals. We divide $480$ minutes of simulation time into $16$ intervals with duration $30$ minutes. The fresh arrival rate is assumed to be piece-wise constants within each interval, but it varies from interval to interval. The fresh arrival pattern is shown in Figure~\ref{fig:arrivalpattern}. This arrival pattern mimics the situation in reality, where there is a morning peak hour and an afternoon peak hour.

We omit the figure for $\mathbf{Z}(t)$, since they are similar to the graph in Figure \ref{fig:validation}. The results for $e_{RD}$ and $e_{RC}$ are shown in Table~\ref{table:scenarios2}.

\begin{table}[ht]
  \begin{center}
  \begin{tabular}{@{}ccc@{}}
    \toprule
    $\hat{\rho}$  & $e_{RD}$ & $e_{RC}$  \\
   	\midrule
  	1.01 &  58.2\% & 1.7\% \\
    1.05 & 32.2\% & 1.4\% \\
    1.1 &  12.8\% & 0.9\% \\
    1.2 &  2.6\% & 0.5\% \\
    1.3 &  1.3\% & 0.4\% \\
    1.4 &  1.1\% & 0.3\% \\
    1.5 &  1.3\% & 0.7\% \\
    \bottomrule
  \end{tabular}
  \caption{Approximation errors of different values of $\hat{\rho}$ in multiple intervals, $\mu=1/4, p=0.5, q=0.1, \theta = 1/2,	\delta_{RD}=0.05, \delta_{RC} = 0.01$.}
  \label{table:scenarios2}
  \end{center}
\end{table}
Similar to Table \ref{table:scenarios1}, one can see from Table~\ref{table:scenarios2} that the fluid model gives close approximations for the number of reconnects for all values of $\hat{\rho}$, and the approximations for the number of redials gets more accurate when $\hat{\rho} > 1.1$.

\section{Erlang A approximation}  \label{erlanga}
The fluid model gives first order approximations for $Z_{Q}(t)$, $Z_{RD}(t)$ and $Z_{RC}(t)$. Based on them, we can approximate the expected total arrival rate and expected number of customers in the queue for any time $t$, from which the expected waiting time can be obtained. However, this is not the eventual goal, since it gives no information about the waiting time distributions of a random customer, which is one of the most used call center performance indicators in call centers. Therefore, to this end, we will apply the Erlang A formula to approximate the waiting time distribution. We assume $\Lambda(t)$ to be the arrival rate of the Erlang A model, which can be obtained via Equation \eqref{eq:totalarrival}.

The reason to use the Erlang A model is intuitively clear, since the redial and reconnect behaviors have only direct influence on the total arrival rate, it has no direct influence on the service, such as the service durations. Therefore, once the total arrival rate $\Lambda(t)$ is given, $Z_{RD}(t)$ and $Z_{RC}(t)$ become irrelevant to what happens in the queue, thus, we can treat the system as an Erlang A system by ignoring the redial and reconnect orbits. Note that this is only an approximation of the Erlang A system, since the arrival process is generally not Poisson.

The analytical expressions for the waiting time distribution and the expected AP of the Erlang A model are known. We refer to \citet{deslauriers2007} and \citet{roubosthesis} for details about the Erlang A formula and calculation
details.

Next, we validate the Erlang A approximation of the original model. To save space, we only show the performances in the case of multiple intervals. The arrival pattern is the same as shown in Figure \ref{fig:arrivalpattern}. For the given parameters, we compute $\mathbf{z}(t)$ via Equations \eqref{eq:fluidlimitq2}-\eqref{eq:fluidlimitrc2}. After that, $\Lambda(t)$ can be obtained via Equation \eqref{eq:totalarrival}. $\Lambda(t)$ will be the input as the arrival rate of the Erlang A formula, from which the SL and AP can be obtained.
The SL is defined as percentage of customers that are answered within $30$ seconds.

We denote SL$^{sim}$ and SL$^a$ as the SL from simulation and from Erlang A formula, respectively. In this section, the SL is specifically set to be the percentage of customers that waited less than $30$ seconds. The AP from simulation and from Erlang A formula are denoted as AP$^{sim}$ and AP$^{a}$, respectively. The comparisons are shown in Table \ref{table:scenarios4}.
\begin{table}[ht]
  \begin{center}
  \begin{tabular}{@{}ccccc@{}}
    \toprule
    $\hat{\rho}$  & $SL^{sim}$ & $SL^{a}$ & $AP^{sim} $ & $AP^{a}$  \\
   	\midrule
    1.01 & 95.5\% & 96.8\% & 5.7\% & 4.8\%  \\
    1.05 & 90.2\% & 92.1\% & 9.1\% & 8.2\%\\
    1.1 & 79.7\% & 80.6\% & 14.2\% & 14.0\% \\
    1.2 & 54.9\% & 53.4\% & 24.2\% & 24.6\% \\
    1.3 & 36.0\% & 34.2\% & 32.9\% & 33.2\% \\
    1.4 & 28.0\% & 26.2\% & 39.8\% & 40.2\% \\
    1.5 & 24.9\% & 23.6\% & 45.6\% & 45.8\% \\
    \bottomrule
  \end{tabular}
  \caption{Approximation errors in SL and AP of different values of $\hat{\rho}$ in multiple intervals, $\mu=1/4, p=0.5, q=0.1, \theta = 1/2, \delta_{RD}=0.05, \delta_{RC} = 0.01$.}
  \label{table:scenarios4}
  \end{center}
\end{table}

Based on the results in Table \ref{table:scenarios4}, we can see that the Erlang A model offers a close approximation both for the SL and AP in all values of $\hat{\rho}$, with the error less than $3\%$ in SL, and $1.1\%$ in the AP.

One might notice that even though we have large errors in $e_{RD}$ when $\hat{\rho} < 1.1$ in Table \ref{table:scenarios1} and Table \ref{table:scenarios2}, the errors in SL and AP are small in Table \ref{table:scenarios4}. This is caused by the fact that when $\hat{\rho}< 1.1$, the number of redials is small compared to the number of reconnects, thus, errors in number of redials would not influence much in $\Lambda(t)$.

\section{Conclusion}
In this paper, we investigate staffing of call centers with redials and reconnects. We consider call centers that operate under heavy load. The model can be described as a three-dimensional Markov process  $\{\mathbf{Z}(t), t>0\}$, defined in \eqref{def1}.
However, to avoid the complexity of solving the Markov process, we use a fluid model to approximate $\mathbf{Z}(t)$.
We show that the fluid limit is the unique solution of a set of three differential equations. Under the same fluid scaling, we derive the fluid limit of the queueing system in the non-stationary case to mimic the real situation in call centers, as the parameters can change before the system reaches stationarity. We also performed simulation experiments to assess the accuracy of the approximations.
To apply the results to real call center applications, we take a further step by calculating the expected total arrival rate, and use this as an input to the Erlang A formula to calculate the SL and AP. Simulation results show that our approximation of the SL is accurate with error less than $2\%$ in all scenarios, and approximation of AP has errors less than $1\%$ when $\hat{\rho} \leq 1.05$ and less than $0.5\%$ when $\hat{\rho} > 1.05$.\\
\ \\
The results suggest a number of topics for further research.
First, the current paper is focused on the derivation and usage of fluid limits for staffing problems of large call centers featuring both redials and reconnects, with load per server greater than $1$. As a next step, it is interesting to supplement the results presented here with the development of staffing methods for the case where the load is strictly less than $1$. To this end, the results of the
present paper and the results for staffing large call centers without redials/reconnects \citet{borst2004, roubosthesis, sze1984} will serve as a good starting point.
Second, with the presence of the redial and reconnect behaviors, it would be interesting to explicitly quantify the reduction of
staffing costs while still meeting the target SL by more efficient planning of call center agents.
Third, besides the influences in call centers staffing, the analysis of reconnect and redial behaviors can also offer insight to call center management. For example, by looking at the reconnect probability of each agent, managers can have some overview information of the quality of service offered by each agent. Furthermore, often the agents have some control on the holding time of each call, and by looking at the correlation between the reconnect probability and the holding time of each call, manager may find the ``right" amount of holding time of each call, such that the holding time and the quality of service is well balanced.

\section*{Appendices}
\subsection*{Appendix A.1: Notations}
~\\
Dividing by $n$ on both sides of Equations \eqref{eq:fluidqn}-\eqref{eq:fluidrcn}, we have
\begin{align}
  \bar{Z}_{Q}^{\left( n\right)} \left( t \right) & = \bar{Z}_Q^{\left( n \right) } \left( 0 \right) +  {G}_{Q}^{\left( n \right)} \left( \mathbf{\bar{Z}}^{\left( n \right)} \right) \left( t \right) + \int_{0}^{t} H_{Q} \left( \mathbf{\bar{Z}}^{\left( n \right)} \right) \left( u \right) du, \label{eq:fluidqndn}\\
  \bar{Z}_{RD}^{\left( n\right)} \left( t \right) & = \bar{Z}_{RD}^{\left( n \right) } \left( 0 \right) + {G}_{RD}^{\left( n \right)}  \left( \mathbf{\bar{Z}}^{\left( n \right)} \right) \left( t\right) + \int_{0}^{t} H_{RD} \left( \mathbf{\bar{Z}}^{\left( n \right)} \right) \left( u\right) du,\label{eq:fluidrdndn} \\
  \bar{Z}_{RC}^{ \left( n \right)} \left( t \right) & = \bar{Z}_{RC}^{\left( n \right) } \left( 0 \right) + {G}_{RC}^{\left( n \right)}\left( \mathbf{\bar{Z}}^{\left( n \right)} \right) \left( t\right) + \int_{0}^{t}H_{RC} \left( \mathbf{\bar{Z}}^{\left( n \right)} \right) \left( u\right) du, \label{eq:fluidrcndn}
\end{align}
where
\begin{align}
{G}_{Q}^{\left( n \right)} \left( \bar{\bf{Z}}^{\left( n \right)} \right) \left( t \right) &: = \Biggl( \frac{{\Pi}_{\lambda n}^{\left( n \right)}\left( t \right)}{n}  - \lambda t \Biggr) - \Biggl( \bar{D}_s^{(n)}(t) -  \int_{0}^{t} \mu \min \{s, \bar{Z}_{Q}^{(n)}(u)\} du \Biggr) \nn \\
  & - \Biggl( \bar{D}_a^{(n)}(t) - \int_0^t \theta \left( \bar{Z}_{Q}^{\left( n \right)} \left( u \right) - s\right)^+ du \Biggr) \nn \\
& + \Biggl( \bar{D}_{RD}^{(n)}\left( t \right) - \int_0^t \delta_{RD} \bar{Z}_{RD}^{\left( n \right)}\left( u \right) du \Biggr) \nn \\
&  + \Biggl( \bar{D}_{RC}^{(n)} \left( t \right) -  \int_0^t \delta_{RC} \bar{Z}_{RC}^{\left( n \right)}\left( u \right) du \Biggr), \label{eq:defgq} \\
G_{RD}^{\left( n \right)} \left( \bar{\bf{Z}}^{\left( n \right)} \right) \left( t \right) & :=  \Biggl( \sum_{j=1}^{n \bar{D}_{a}^{\left( n \right)} \left( t\right)} B_j\left( p\right)/ n - \int_0^t p\theta \left( \bar{Z}_{Q}^{\left( n \right)} \left( u \right) - s\right)^+ du \Biggr) \nn \\
&- \Biggl( \bar{D}_{RD}^{\left( n \right)} \left( t\right)
- \int_0^t \delta_{RD} \bar{Z}_{RD}^{\left( n \right)}\left(u\right) du \Biggr), \label{eq:defgrd} \\
 {G}_{RC}^{\left( n \right)} \left( \bar{\bf{Z}}^{\left( n \right)} \right) \left( t \right) & := \Biggl( \sum_{j=1}^{n \bar{D}_{s}^{\left( n \right)} \left( t\right)} B_j\left( q\right) / n - \int_0^t q \mu \min \{s, \bar{Z}_{Q}^{(n)}(u) \} du \Biggr) \nn \\
   & - \Biggl( \bar{D}_{RC}^{(n)} \left( t \right) -  \int_0^t \delta_{RC} \bar{Z}_{RC}^{\left( n \right)}\left( u \right) du \Biggr), \label{eq:defgrc}
\end{align}
and
\begin{align}
\bar{D}_s^{(n)}(t) & = \Pi_1 \left( n \int_{0}^{t} \mu \min \{s, \bar{Z}_{Q}^{(n)}(u)\} du \right) /n , \label{eq:ds} \\
\bar{D}_a^{(n)}(t) & = \Pi_2 \left( n \int_0^t \theta \left( \bar{Z}_{Q}^{\left( n \right)} \left( u \right) - s\right)^+ du \right) /n, \label{eq:da} \\
\bar{D}_{RD}^{(n)}(t) & = \Pi_3 \left( n \int_0^t \delta_{RD} \bar{Z}_{RD}^{\left( n \right)}\left( u \right) du \right) /n, \label{eq:drd} \\
\bar{D}_{RC}^{(n)}(t) & = \Pi_4 \left( n \int_0^t \delta_{RC} \bar{Z}_{RC}^{\left( n \right)}\left( u \right) du \right) /n, \label{eq:drc}
\end{align}
and
\begin{align}
  \int_{0}^{t}H_{Q} \left( \bar{\bf{Z}}^{\left( n \right)} \right) \left( u \right) du &: = \int_{0}^{t}\lambda +  \delta_{RD} \bar{Z}_{RD}^{\left( n \right)} \left( u \right) + \delta_{RC} \bar{Z}_{RC}^{\left( n \right)} \left( u \right) \nonumber \\
  & - \mu \min \{s, \bar{Z}_{Q}^{\left( n \right)} \left( u \right)\} - \theta\left( \bar{Z}_{Q}^{\left( n \right)} \left( u \right) - s \right)^+ du, \nn \\
  \int_{0}^{t}H_{RD} \left( \bar{\bf{Z}}^{\left( n \right)} \right) \left( u\right) du & :=  \int_{0}^{t} p \theta\left( \bar{Z}_{Q}^{\left( n \right)} \left( u \right) - s \right)^+ - \delta_{RD} \bar{Z}_{RD}^{\left( n \right)} \left( u \right) du,  \nn \\
  \int_{0}^{t}H_{RC} \left( \bar{\bf{Z}}^{\left( n \right)} \right) \left( u\right) du & := \int_{0}^{t} q \mu \min \{s, \bar{Z}_{Q}^{\left( n \right)} \left( u \right)\} - \delta_{RC} \bar{Z}_{RC}^{\left( n \right)} \left( u \right) du. \nn
\end{align}

For the convenience of notation, we rewrite Equations \eqref{eq:fluidqndn}-\eqref{eq:fluidrcndn} in the vector form
\begin{align}
\mathbf{\bar{Z}}^{\left( n\right)} \left( t \right) & = \mathbf{\bar{Z}}^{\left( n \right) } \left( 0 \right)  +   \mathbf{G}^{\left( n \right)} \left(\bar{\bf{Z}}^{\left( n \right)} \right)\left(t\right) + \int_{0}^{t}\mathbf{H} \left(\bar{\bf{Z}}^{\left( n \right)} \right)\left(u\right)   du, \label{eq:fluidvector}
\end{align}
where
\begin{align*}
& \mathbf{{G}}^{\left( n \right)} \left(\bar{\bf{Z}}^{\left( n \right)} \right)\left(t\right) := \left( {G}_Q^{\left( n \right)}\left(\bar{\bf{Z}}^{\left( n \right)} \right) \left(t\right),  {G}_{RD}^{\left( n \right)}\left(\bar{\bf{Z}}^{\left( n \right)} \right) \left(t\right), {G}_{RC}^{\left( n \right)}\left(\bar{\bf{Z}}^{\left( n \right)} \right) \left(t\right) \right)^T, \\
 & \mathbf{H}\left(\bar{\bf{Z}}^{\left( n \right)} \right)\left(u\right) := \left( H_Q\left(\bar{\bf{Z}}^{\left( n \right)} \right) \left(u\right),  H_{RD}\left(\bar{\bf{Z}}^{\left( n \right)} \right)\left(u\right), H_{RC}\left(\bar{\bf{Z}}^{\left( n \right)}\right) \left(u\right) \right)^T.
\end{align*}

\bigskip

\subsection*{Appendix A.2: Proof of Lemma~\ref{lemma}}
\begin{proof}
In order to show that $\displaystyle \{\mathbf{\bar{Z}}^{\left( n \right)}( \cdot)\}_{n=1}^{\infty}$ is relatively compact with continuous limits, it is sufficient to show the following two properties (see Corollary 7.4 and Theorem 10.2 of \citet{ethier}).
\begin{enumerate}
\item Compact Containment:
for any $T \geq 0, \epsilon > 0$, there exists a compact set $\Gamma_T \subset  \mathbb{R}^3$ such that
\begin{align*}
  \p\left( \mathbf{\bar{Z}}^{(n)}(t) \in \Gamma_T, \ t \in [0,T] \right) \rightarrow 1, \hspace{3mm} \text{as } n \rightarrow \infty;
\end{align*}
\item Oscillation Control:
for any $\epsilon>0$, and $T \geq 0$, there exists a $\delta > 0$, such that
\begin{align}
\limsup_{n \to \infty} \p\left( \omega\left( \mathbf{\bar{Z}}^{(n)}, \delta, T \right) \geq \epsilon \right) \leq \epsilon, \label{eq:oscillation_control}
\end{align}
where
\begin{align*}
\omega(\mathbf{x}, \delta, T):= \sup_{\substack{\nu, t \in [0, T] \\ |s-\nu|<\delta}} \max_{j \in J}|x_j(t) - x_j(\nu) |,
\end{align*}
and $J := \{Q, RD, RC\}$.
\end{enumerate}

{\raggedleft Proof of Compact Containment property:} \\
The following trivial upper bound holds for the total number of customers in the system (only arrivals are taken into account and no departures): for $t \in [0, T]$,
\begin{align*}
  \bar{Z}^{(n)}_{Q}(t)+ \bar{Z}^{(n)}_{RD}(t)+ \bar{Z}^{(n)}_{RC}(t) \leq  \bar{Z}^{(n)}_{Q}(0)+ \bar{Z}^{(n)}_{RD}(0)+ \bar{Z}^{(n)}_{RC}(0) + \Pi_{\lambda n}^{(n)}(T) /n.
\end{align*}

Since $\Pi_{\lambda n}^{(n)}(\cdot) $ is a Poisson process of rate $\lambda n$, by the Law of Large Numbers(LLN), we have
\begin{align*}
\Pi_{\lambda n}^{(n)}(T)/n \xrightarrow{d} \lambda T \hspace{3mm} \text{as } n \rightarrow \infty.
\end{align*}

By the assumption of Theorem~\ref{theo:theo1}, we have
\begin{align*}
 \bar{Z}_Q^{\left( n\right)}(0) + \bar{Z}_{RD}^{\left( n\right)}(0) + \bar{Z}_{RC}^{\left( n\right)}(0) \xrightarrow{d} z_{Q}(0) + z_{RD}(0) + z_{RC}(0) .
\end{align*}

Hence
\begin{align*}
  \p( \mathbf{\bar{Z}}^{(n)}(t) \in \Gamma_T, \ t \in [0,T] ) \rightarrow 1 \hspace{3mm} \text{as } n \rightarrow \infty,
\end{align*}
where $\Gamma_T = \{ (x_1, x_2, x_3) \mid x_1+x_2+x_3 \leq z_Q(0)+ z_{RD}(0)+ z_{RC}(0) + \lambda T + 1, \  x_1, x_2, x_3 \geq 0 \}$, and the compact containment property indeed holds.
\bigskip

{\raggedleft Proof of Oscillation Control property:} \\
It follows from Equations \eqref{eq:fluidqn}-\eqref{eq:fluidrcn} that, for all $\nu, t \geq 0$,
\begin{align*} 
 | \bar{Z}_{Q}^{(n)}(t) - \bar{Z}_{Q}^{(n)}(\nu) | &\leq  |\Pi_{\lambda n}^{(n)}(t) - \Pi_{\lambda n}^{(n)}(\nu)| / n  + \sum_{j \in \{ s, a, RD, RC \}} |\bar{D}_j^{(n)}(t) - \bar{D}_j^{(n)}(\nu)|, \nn \\
| \bar{Z}_{RD}^{(n)}(t) - \bar{Z}_{Q}^{(n)}(\nu) | &\leq  \sum_{j \in \{ a, RD\}} |\bar{D}_j^{(n)}(t) - \bar{D}_j^{(n)}(\nu)|, \nn \\
| \bar{Z}_{RC}^{(n)}(t) - \bar{Z}_{Q}^{(n)}(\nu) | &\leq  \sum_{j \in \{ s, RC \}} |\bar{D}_j^{(n)}(t) - \bar{D}_j^{(n)}(\nu)|,
\end{align*}
where the processes $\bar{D}_j^{(n)}(\cdot)$ are defined by \eqref{eq:da}-\eqref{eq:drc}.

Also, from the Compact Containment property, we know that there exists a finite constant $V$ such that
\begin{align*}
\p \left( \underbrace{\bar{Z}_{Q}^{(n)}(u), \bar{Z}_{RD}^{(n)}(u), \bar{Z}_{RC}^{(n)}(u) \leq V,  \ u \in [0,T] }_{\displaystyle {=: \Omega_n}} \right) \to 1 \hspace{3mm} \text{as } n \rightarrow \infty.
\end{align*}

On the event $\Omega_n$, the following inequalities hold for all $\nu, t \in [0,T]$ such that $|t - \nu| \leq \delta$:
\begin{equation*} 
  \begin{aligned}
  & \int_{\nu}^{t}\mu \min\{s, \bar{Z}_{Q}^{(n)}(u) \}du \leq \mu s \delta = : c_1 (\delta), \\
  & \int_{\nu}^{t}\theta \left( \bar{Z}_{Q}^{(n)}(u) - s\right)^+  du \leq \theta V \delta = : c_2 (\delta), \\
  & \int_{\nu}^{t}\delta_{RD}\bar{Z}_{RD}^{(n)}(u)du \leq \delta_{RD}V \delta = : c_3 (\delta), \\
  & \int_{\nu}^{t}\delta_{RC}\bar{Z}_{RC}^{(n)}(u)du \leq \delta_{RC}V \delta = : c_4 (\delta).
  \end{aligned}
\end{equation*}

Employing formulas \eqref{eq:ds}-\eqref{eq:drd}, we then get
\begin{align*}
\p\left( \omega\left( \mathbf{\bar{Z}}^{(n)}, \delta, T \right) \geq \epsilon \right) & \leq \p\left( \omega\left( \Pi_{\lambda n}^{(n)}(\cdot)/n, \delta, T \right) \geq \epsilon /5 \right) \nn \\
&+ \sum_{j=1}^4 \p\left( \omega\left( \Pi_j(n\cdot)/n, c_j(\delta), c_j(\delta )T \right) \geq \epsilon/5  \right) \nn \\
& \to \lambda \delta + \sum_{j=1}^4 c_j(\delta) = \delta (\lambda + \mu s + \delta_{RD} V + \delta_{RC} V),
\end{align*}
where the convergence holds by the LLN for the Poisson processes $\Pi_{\lambda n}^{(n)}(\cdot)/n$, $\Pi_j(n\cdot)/n$ and by the continuity of the moduli of continuity $\omega(x(\cdot), \delta, T)$, $\omega(x(\cdot), c_j(\delta), c_j(\delta)T)$ with respect to $x(\cdot)$.

Hence, the oscillation control property \eqref{eq:oscillation_control} indeed holds with $\delta  = \epsilon / (\lambda + \mu s + \delta_{RD} V + \delta_{RC} V)$.

\end{proof}
\bigskip

\subsection*{Appendix A.3: Proof of Theorem~\ref{theo:theo1} }
\begin{proof}
In Lemma~\ref{lemma}, we have shown that the sequence $\{\bar{\mathbf{Z}}^{n}(\cdot)\}_{n=1}^{\infty}$ is relatively compact with continuous limits, that is, from any subsequence $\{\bar{\mathbf{Z}}^{n_k}(\cdot)\}_{k=1}^{\infty}$, we can extract another subsequence $\{\bar{\mathbf{Z}}^{n_{k_l}}(\cdot)\}_{l=1}^{\infty}$ that converges weakly in $D([0, \infty),  \mathbb{R}^3)$, say to a continuous process $z^{*}(t)$. We then call $z^{*}(t)$ a particular limit of the original sequence $\{\bar{\mathbf{Z}}^{n}(\cdot)\}_{n=1}^{\infty}$.

Consider an arbitrary particular limit $\mathbf{z}^*(\cdot)$ along a subsequence $\{\bar{\mathbf{Z}}^{n_k}(\cdot)\}_{k=1}^{\infty}$. If we can show that $\mathbf{z}^*(\cdot)$ satisfies Equations \eqref{eq:fluidlimitq}-\eqref{eq:fluidlimitrc}, and Equations \eqref{eq:fluidlimitq}-\eqref{eq:fluidlimitrc} have a unique solution, then, due to the arbitrariness of $\mathbf{z}^*(\cdot)$, there must be a unique fluid limit defined by Equations \eqref{eq:fluidlimitq}-\eqref{eq:fluidlimitrc}.

We have
\begin{equation} \label{eq:conv1}
\bar{\mathbf{Z}}^{n_k}(\cdot) - \bar{\mathbf{Z}}^{n_k}(0) - \int_0^\cdot \mathbf{H} \left( \bar{\mathbf{Z}}^{n_k} \right) du = \mathbf{G}^{n_k}(\cdot).
\end{equation}

On one hand, since $\bar{\mathbf{Z}}^{n_k}(\cdot) \xrightarrow{d} \mathbf{z}^*(\cdot)$ as $k \to \infty$ and the limit $\mathbf{z}^*(\cdot)$ is continuous,
\[
\bar{\mathbf{Z}}^{n_k}(\cdot) - \bar{\mathbf{Z}}^{n_k}(0) - \int_0^\cdot \mathbf{H} \left( \bar{\mathbf{Z}}^{n_k} \right) du \xrightarrow{d} \mathbf{z}^\ast(\cdot) - \mathbf{z}(0) - \int_0^\cdot \mathbf{H} \left( \mathbf{z}^\ast \right) du.
\]

On the other hand, below we show that $\mathbf{G}^{n_k}(\cdot) \xrightarrow{d} 0$, and then \eqref{eq:conv1} implies that
\[
\bar{\mathbf{Z}}^{n_k}(\cdot) - \bar{\mathbf{Z}}^{n_k}(0) - \int_0^\cdot \mathbf{H} \left( \bar{\mathbf{Z}}^{n_k} \right) du \xrightarrow{d} 0.
\]

As we combine the last two displays together, it follows that the particular limit $\mathbf{z}^\ast$ a.s. satisfies Equations \eqref{eq:fluidlimitq}-\eqref{eq:fluidlimitrc}. Also, the mapping $\mathbf{H}$ is Lipschitz continuous and then, by Lemma 1 in~\citet{reed2004diffusion}, Equations \eqref{eq:fluidlimitq}-\eqref{eq:fluidlimitrc} have a unique solution. Hence, all particular fluid limits are the same, namely they coincide with the unique solution to \eqref{eq:fluidlimitq}-\eqref{eq:fluidlimitrc}.

\bigskip
It is left to show that $\mathbf{G}^{n_k}(\cdot) \xrightarrow{d} 0$.

By the LLN,
\[
\Pi_1(n \cdot) / n - \cdot \xrightarrow{d} 0 \ \text{in} \ D([0,\infty), \mathbb{R}),
\]
and also, since $\bar{\mathbf{Z}}^{n_k}(\cdot) \xrightarrow{d} \mathbf{z}^*(\cdot)$ and $\mathbf{z}^\ast$ is continuous,
\[
 \int_0^\cdot \mu \min\{s, \bar{Z}_{Q}^{(n_k)}(u) \}du \xrightarrow{d} \int_0^\cdot \mu \min\{s, \mathbf{z}^\ast (u) \}du  \ \text{in} \ D([0,\infty), \mathbb{R}).
\]
Then, by \eqref{eq:ds} and the Random time change Theorem in~\citet{billingsley2009convergence},
\[
\bar{D}_s^{(n_k)}(t) -  \int_{0}^{t} \mu \min \{s, \bar{Z}_{Q}^{(n_k)}(u)\} du \xrightarrow{d} 0 \ \text{in} \ D([0,\infty), \mathbb{R}).
\]

By the same argument, one can show that the other terms in $G_Q^{(n_k)}(\cdot)$ converge to $0$, and that $G_{RC}^{(n_k)}(\cdot)$, $G_{RD}^{(n_k)}(\cdot)$ converge to $0$, too. Hence, the proof of Theorem~\ref{theo:theo1} is finished.
\end{proof}
\newpage
\bibliographystyle{plainnat}

\end{document}